\documentclass[12pt]{article}
\usepackage{a4}
\usepackage{latexsym}
\usepackage{amsmath}
\usepackage{amssymb}
\usepackage{amscd}
\usepackage[all]{xy}
\usepackage{amsthm}

\theoremstyle{plain}
    \newtheorem{theorem}                    {Theorem}       [section]
    \newtheorem{lemma}      [theorem]       {Lemma}
    \newtheorem{corollary}  [theorem]       {Corollary}
    \newtheorem{proposition}[theorem]       {Proposition}
    
    \newtheorem{conjecture} [theorem]       {Conjecture}

\newtheorem{remark}[theorem]{Remark}

\DeclareFontFamily{U}{wncy}{}
    \DeclareFontShape{U}{wncy}{m}{n}{<->wncyr10}{}
    \DeclareSymbolFont{mcy}{U}{wncy}{m}{n}
    \DeclareMathSymbol{\Sh}{\mathord}{mcy}{"58} 

\newcommand{\chr}{\operatorname{char}}
\newcommand{\Hom}{\operatorname{Hom}}
\newcommand{\Gal}{\operatorname{Gal}}
\newcommand{\Pic}{\operatorname{Pic}}

\newcommand{\Tor}{{\operatorname{Tor}}}

\newcommand{\Spec}{\operatorname{Spec}}

\newcommand{\im}{\operatorname{im}}

\newcommand{\coker}{\operatorname{coker}}
\newcommand{\id}{\operatorname{id}}
\newcommand{\Br}{\operatorname{Br}}
\newcommand{\CH}{\operatorname{CH}}
\renewcommand{\lim}{\operatornamewithlimits{lim}}
\newcommand{\colim}{\operatornamewithlimits{colim}}

\newcommand{\NS}{\operatorname{NS}}
\newcommand{\f}{{\mathcal F}}
\newcommand{\g}{{\mathcal G}}

\newcommand{\Z}{{{\mathbb Z}}}

\newcommand{\ZZ}{{{\textstyle\mathbb Z'}}}
\newcommand{\Q}{{{\mathbb Q}}}

\newcommand{\F}{{{\mathbb F}}}

\newcommand{\T}{{{\mathcal T}}}

\newcommand{\X}{{\mathcal X}}
\newcommand{\Y}{{\mathcal Y}}

\newcommand{\M}{{\mathcal M}}

\newcommand{\et}{{\text{\rm et}}}

\newcommand{\sep}{{\text{\rm s}}}

\newcommand{\Alb}{\operatorname{Alb}}
\newcommand{\Ext}{\operatorname{Ext}}
\newcommand{\Sel}{\operatorname{Sel}}

\newcommand{\corank}{\operatorname{corank}}
\renewcommand{\O}{{\cal O}}
\newcommand{\red}{\text{red}}

\date{}

\title{Tate's conjecture and the Tate-Shafarevich group over global
function fields.}

\author{Thomas H. Geisser}

\begin{document}
\maketitle

\begin{abstract}
Let $\X$ be a regular variety, flat and proper over a complete regular 
curve over a finite field, such that the generic fiber
$X$ is smooth and geometrically connected. 
We prove that the Brauer group of $\X$ is finite if
and only Tate's conjecture for divisors on $X$ holds and the
Tate-Shafarevich group of the Albanese variety of $X$ is finite,
generalizing a theorem of Artin and Grothendieck for surfaces
to arbitrary relative dimension.
We also give a formula relating the orders of the group under
the assumption that they are finite, generalizing the formula
given for a surface in \cite{ichsurface}.
\end{abstract}

\section{Introduction}
Let $C$ be a smooth and proper 
curve over a finite field of characteristic $p$ with function field $K$. 
Let $\X$ be a regular scheme and $\X\to C$ a proper flat map 
such that $X=\X\times_CK$ is smooth with geometrically connected fibers over $K$.
Let $K_v$ be the quotient field of the henselization of $C$ at a point $v$, and 
$X_v=X\times_KK_v$. 

If $\X$ is a surface, then it is a classical result of Artin
and Grothendieck \cite[\S 4]{brauerIII} that the Brauer group of $\X$ is finite 
if and only if the Tate-Shafarevich group of the Jacobian $A$ of $X$ is finite.
Moreover, if they are finite, then their orders are related by the formula \cite{ichsurface}
\begin{equation}\label{lformul}
|\Br(\X)|\alpha^2  \delta^2= |\Sh(A)|\prod_{v\in V} \alpha_v\delta_v,
\end{equation}
where $\delta$ and $\delta_v$ are the indices of $X$ and $X_v$,
respectively, and $\alpha$ and $\alpha_v$ are the order of 
the cokernel of the inclusion $\Pic^0(X)\to H^0(K,\Pic_X^0)$
and $\Pic^0(X_v)\to H^0(K_v,\Pic_X^0)$, respectively. 
The purpose of this paper to generalize these results 
to arbitrary relative dimension. 

\begin{theorem}\label{thm1}
The Brauer group of $\X$ is finite if and only if 
the Tate-Shafarevich group of the Albanese variety 
of $X$ is finite and Tate's conjecture for divisors holds for $X$, i.e., for all
$l\not=p$, the cycle map
$$c_l: \Pic(X)\otimes\Z_l \to H^2_\et(X^s,\Z_l(1))^{\Gal(K)}$$
has torsion cokernel, where $X^s$ is the base extension to the separable
closure.
\end{theorem}

More precisely, we have 
$\corank(\coker c_l), \corank(\Sh(\Alb_X))\leq \corank(\Br(\X))$,
and $\corank(\Br(\X))\leq \corank(\coker c_l)+ \corank( \Sh(\Alb_X))$.
If $\X$ is a surface, then Tate's conjecture for divisors on $X$  
is trivial and we recover the classical result of Artin-Grothendieck. 

In order to generalize the formula \eqref{lformul}, we need to introduce 
other invariants which vanish for a relative curve. Consider the maps
\begin{align*}
\xi_n^i : H^2(K,H^i_\et(X^s,\Q/\Z[{\textstyle\frac{1}{p}}](n)))&\longrightarrow
\bigoplus_v H^2(K_v,H^{i}_\et(X^s,\Q/\Z[{\textstyle\frac{1}{p}}](n))), \\
l_n^i: H^i_\et(X,\Z(n)) & \longrightarrow \prod_v   H^i_\et(X_v,\Z(n)),
\end{align*}
where $H^i_\et(X,\Z(n))$ is the etale hypercohomology of Bloch's cycle complex.  
The maps $\xi_n^i$ were studied by Jannsen \cite{jannsenggq}, and the 
argument of loc.cit. Thm.\ 3 shows that they have finite kernel and cokernel 
if $i\not=2n-2$. 

Let $\beta, \beta_v$ be the orders of the prime to $p$-part of the cokernels of 
$$H^{2d}_\et(X,\Z(d))\stackrel{\rho_d}{\longrightarrow} H^{2d}_\et(X^s,\Z(d))^{\Gal(K)}, \quad
H^{2d}_\et(X_v,\Z(d))\stackrel{}{\longrightarrow} H^{2d}_\et(X_v^\sep,\Z(d))^{G_v}$$
and let $\delta', \delta'_v$ be the order of 
prime to $p$-part of the cokernels of the maps 
$$CH^d(X^\sep)^G\stackrel{\deg}{\longrightarrow} \Z, \quad 
CH^d(X_v^\sep)^{G_v}\stackrel{\deg}{\longrightarrow} \Z, $$ respectively. 
Then $\alpha, \delta'$, and $\prod_v \delta'_v$ are finite, and 
we use a result of Saito-Sato to show that $\prod_v\beta_v$ 
is finite. Let $c$ the product of the maps $c_l$
for all $l\not=p $. 

\begin{theorem} \label{order}
Let $X$ be smooth and projective over a global field
of characteristic $p$. If $X$ admits a regular proper model, 
$\Sh(\Alb_X)$ is finite, Tate's conjecture for divisors holds
on $X$ and $\ker \xi_d^{2d-2}$ is finite, 
then up to a power of $p$ we have
$$|\ker l_d^{2d+1}|\cdot |\coker c|\cdot \alpha\beta\delta'=
|\Sh(\Alb_X)|\cdot |\ker \xi_d^{2d-2}|\prod_v \beta_v\delta'_v.$$
\end{theorem}

Note that the terms appearing in the theorem only depend on $X$. If 
$\ker \xi_d^{2d-2}$, or equivalently $\coker \rho_d$, is not finite, 
then the theorem still holds if we replace the corresponding terms by 
the order of their finite 
quotients by a common subgroup. If $X$ is a curve, then the formula 
reduces to \eqref{lformul} because $\ker \xi_1^0=0$ and 
$\ker l_1^3\cong \Br(\X)$.

The same argument gives the following theorem which relates the above
mentioned groups,
and implies Theorem \ref{thm1}. Let $\Alb_X$ be the Albanese variety of $X$, 
and for an abelian group $A$ let $A^*$ and 
$TA$ be the Pontrjagin dual and Tate-module, respectively. 

\begin{theorem}\label{maines}
Modulo the Serre subcategory of finite groups and $p$-power torsion groups, 
we have an exact sequence of torsion groups
\begin{multline*}
 0\to\coker \rho_d\to  \ker \xi_d^{2d-2} \to \ker l_d^{2d+1}\to\Sh(\Alb_X)\\
\to (\coker c)^* \to (T\Br(\X))^*\to (T\Sh(\Pic_X^0))^*\to 0.
\end{multline*}
\end{theorem}


During the proof we obtain the following result on the maps $\xi_d^{2d-2}$ and $l_d^{2d+1}$.

\begin{proposition}\label{thm3}\label{secondes}
Assume that $\Br(\X)$ is finite.

1)  There is a complex 
\begin{multline*}
0 \to \coker\rho_d \to 
H^2(K,H^{2d-2}_\et(X^s,\Q/\Z(d))) \stackrel{\xi_d^{2d-2}}{\longrightarrow}\\
\bigoplus_v H^2(K_v,H^{2d-2}_\et(X^s,\Q/\Z(d)))\to \Hom(\NS(X),\Q/\Z)\to 0
\end{multline*}
which is exact up to finite groups and $p$-groups.

2) There is a  complex 
\begin{multline*}
H^{2d+1}_\et(\X,\Z(d)) \to 
H^{2d+1}_\et(X,\Z(d))\stackrel{l_d^{2d+1}}{\longrightarrow}\\
\bigoplus_v H^{2d+1}_\et(X_v,\Z(d))\to \Hom(\Pic(X),\Q/\Z)\to 0
\end{multline*}
which is exact up to $p$-groups.
\end{proposition}

\medskip

The idea of the proof of Theorem \ref{maines} is to work with one-dimensional
cycles and to use a theorem of Saito-Sato \cite{saisa}.
For our purposes,
it is necessary to give a slight improvement of their main theorem,
which we are able to prove using results of Gabber \cite{gabberthm} 
and Kerz-Saito \cite{KS}:
Let $f:\Y\to \Spec R$ be of finite type over the spectrum of an excellent 
henselian discrete valuation ring $R$ with residue field $k$.
The following conjecture and theorem on Kato-homology
were stated and proven in \cite{saisa}
for $\Y$ projective over $R$ such that the reduced special fiber
is a strict normal crossing scheme. 

\begin{conjecture}\cite[Conj. 2.11, 2.12]{saisa}
Let $\Y$ be a regular scheme, flat and proper over $\Spec R$,
and $l$ a prime invertible on $R$.

1) If $k$ is separably closed, then $KH_a(\Y,\Q_l/\Z_l)=0$ for all $a$.

2) If $k$ is finite, then 
$$KH_a(\Y,\Q_l/\Z_l)=
\begin{cases}
0  &\text{for } a\not=1;\\
(\Q_l/\Z_l)^{I} &\text{for }a=1,
\end{cases} $$
where $I$ is the number of irreducible components of the special fiber.
\end{conjecture}

\begin{theorem}\cite[Thm. 2.13]{saisa}
The conjecture holds for $a\leq 3$.
\end{theorem}

We thank T. Szamuely for helpful comments on 
a preliminary version of this paper, and K. Sato for discussions on his work.

\medskip 
{\bf In the body of this paper we invert the characteristic $p$ of the 
base field, i.e.,
we work with $\Z'=\Z[\frac{1}{p}]$-coefficients}.
For a compact or discrete abelian group $A$ we define 
the Pontrjagin dual to be $A^*=\Hom_{cont}(A, \Q/\Z[\frac{1}{p}])$, 
the completion to be $A^\wedge= \lim_{p\not|m} A/m$,
the torsion to be $\Tor A= \colim_{p\not|m} {}_mA$, 
and the Tate-module to be $TA=\lim_{p\not|m}{}_mA$.
Note that the results for $A$ and $A\otimes \Z'$  agree, 
so that we sometimes omit the $-\otimes\Z'$.

\section{Etale motivic cohomology}
Let $f:\Y\to B$ be separated and of finite type over the spectrum
$B$ of a Dedekind ring of exponential characteristic $p$.
Consider Bloch's cycle complex $\Z^c(w)$ 
of cycles of relative dimension $w$ with  
the etale sheaf $z_w(-,-i-2w)$ in (cohomological) degree $i$.
For 
an abelian group $A$, we define and etale motivic Borel-Moore homology 
$H_{i}^c(\Y,A(w))$ to be $H^{-i}(\Y,A\otimes \Z^c(w))$. 
For finite maps $g:\Y\to \X$, the
direct image functor is exact, $g_*=Rg_*$. This implies that
proper push-forward of 
algebraic cycles induces covariant functoriality. 
Flat pull-back induces contravariant functoriality with 
the appropriate shift in weight. 
Since the homology of Bloch's complex
agrees with its etale hypercohomology with rational coefficients, 
the natural map 
$$CH_w(\Y,i-2w) \to H_i^c(\Y,\Z(w))  $$ 
is an isomorphism upon tensoring with $\Q$, 
and this vanishes if $i<2w$. In particular, we have
$H_i^c(\Y,\Z(n)) \cong H_{i+1}^c(\Y,\Q/\Z(n))$ for $i<2w-1$.

If $\Y$ is regular and of pure dimension $d+1$, 
we also write $\Z(n)$ for $\Z^c(d+1-n)[-2d-2]$ 
and $H^{i}_\et(\Y,A(n))$ for $H_{2d+2-i}^\et(\Y,A(d+1-n))$. 
If, in addition, $B$ is a regular curve over a perfect field and 
$\Y$ regular, then 
$\Z^c/m(d+1-n)[-2d-2]\cong \Z/m(n)\cong \mu_m^{\otimes n}$ 
for $n\geq 0$ and $m$ invertible on $\Y$, 
because this holds for smooth schemes over fields \cite{ichmarcII}.

By \cite[XVIII, Th. 3.1.4]{SGA4}, the functor $Rf_!$ has a right
adjoint $Rf^!$ on the category of bounded complexes of $\Z/m$-modules
for any compactifiable morphism $f$.

\begin{lemma}\%label{ssff}
If $B$ is a regular curve over a perfect field, 
and $\Y$ is regular, then 
$$Rf^!\mu_m^{\otimes -n} \cong \mu_m^{\otimes d-n}[2d]\cong \Z^c/m(n+1)[-2]$$
for any $m$ invertible on $B$ and $n\leq d$.
\end{lemma}

Taking the limit over etale neighborhoods of a point of $B$, the 
Lemma remains true for $B$ the henselization of a regular curve
over a perfect field.

\proof
Let $s:B\to k$ be the structure map. Then $s$ and $f$ are smooth,
and $ Rs^!\mu_m^{\otimes -n-1}[-2]\cong \mu_m^{\otimes -n}$
on $B$ by \cite[XVIII Thm. 3.2.5]{SGA4}. Hence for the same reason
$$ Rf^!\mu_m^{\otimes -n} \cong R(sf)^!\mu_m^{\otimes -n-1}[-2]
\cong \mu_m^{\otimes d-n}[2d] \cong \Z^c/m(n+1)[-2].$$
\endproof

In particular, we obtain an isomorphism 
$$H_a^c(\Y,\Z/m(1))\cong H_a^\et(\Y,\Z/m(1)):=
H^{2-a}(\Y, Rf^!\Z/m),$$
i.e. for regular schemes, etale motivic Borel-Moore cohomology 
is isomorphic to the etale motivic homology used in \cite{saisa}.

\begin{theorem}\label{descc}
If $T$ is of finite type over a separably closed field $k$ and $w\leq 0$, 
then the natural map 
$$CH_w(T,i-2w)_{\Z'} \to H_{i}^c(T,\ZZ(w))$$
is an isomorphism. 
\end{theorem}

\proof
The proof of \cite[Thm. 3.1]{ichduality} works in this situation, 
but we replicate it for the convenience of the reader. 
Since  $\Z^c(w)$ satisfies the localization property, 
we can apply Thomason's argument \cite[Prop. 2.8]{thomason} 
using induction on the dimension of $T$ 
to reduce to showing that for an artinian local ring $R$, 
essentially of finite type over $k$, the canonical map
$\ZZ^c(w)(\Spec R) \to R\Gamma_\et(\Spec R,\ZZ^c(w))$
is a quasi-isomorphism.
Since $\ZZ^c(w)(U)\cong \ZZ^c(w)(U^{red})$ we can assume that R is reduced, 
in which case it is the spectrum of a field $F$ of finite transcendence 
degree $d$ over $k$. 
We have to show that the canonical map 
$H_i(\ZZ(w)(F))\cong H_i^c(F,\ZZ(w))$ is an 
isomorphism for all $i$. Rationally, Zariski and etale hypercohomology of 
the cycle complex agree. With prime to $p$-coefficients, 
both sides agree for $i\geq d+w $ by the Rost-Voevodsky theorem, and 
for $i<d$ both sides vanish because 
$H_i^c(F,\Z/l(w))\cong H^{2d-i}_\et(F,\Z/l(d-w))$
and the $l$-cohomological dimension of $F$ is $d$. 
\endproof


\begin{proposition}\label{rojtman}
Let $T$ be a smooth projective scheme over a separably closed field $k$.
Then $H^i_\et(T,\ZZ(d))=0$ for $i>2d$, 
$CH_0(T)_{\Z'}\cong H^{2d}_\et(T,\ZZ(d))$, and the albanese map induces
a surjection from the degree zero part
$$ CH_0(T)_{\Z'}^0\to \Alb_T(k)_{\Z'}$$
with uniquely divisible kernel.
\end{proposition}

\proof
The first two statement follow from Theorem \ref{descc}. 
By Rojtman's theorem, the albanese map induces an isomorphism on prime 
to $p$ torsion subgroups \cite[Thm. 4.2]{colthe}.
Hence the final
statement follows by observing that the albanese map is surjective
and that $CH_0(T)^0$ is divisible by all
integers prime to $p$.
\endproof

The following proposition can be proved as \cite[Thm. 1.1]{ichstructure}:

\begin{proposition}\label{strstr}
Let $T$ be smooth and projective over a separably closed field. 
Then $H^i_\et(T,\Z(n))\otimes\Q/\ZZ=0$ for $i\not=2n$. 
\end{proposition}

\subsection{Duality}
We recall some facts on duality from \cite{ichtata}, \cite{ichduality}. 
If $g:Y\to \F$ is separated and of finite type over a finite field
and $m$ prime to the characteristic of $\F$, 
then the adjunction $Rg_! \vdash Rg^!$ gives
\begin{multline*}
 R\Hom_{Ab}(R\Gamma_c(Y,\Z/m(n)),\Z/m[-1]) \cong 
R\Hom_{\F}(Rg_!\Z/m(n),\Z/m)  \\
\cong R\Hom_Y(\Z/m,Rg^!\Z/m(-n)) \cong R\Gamma(Y,Rg^!\Z/m(-n)),
\end{multline*}
where the shift in degree appears because the left derived functor
of coinvariants under $G_{\F}$ is the right derived of the invariant
functor shifted by one. Using the usual definition
$H_i^\et(Y,\Z/m(n)):= H^{-i}_\et(Y,Rg^!\Z/m(-n))$ for etale homology, 
we obtain a perfect pairing:
$$ H^i_c(Y,\Z/m(n)) \times H_{i-1}^\et(Y,\Z/m(n))\to   \Z/m.$$
For finite maps $\iota :Y \to \X$, the adjunction maps
$\iota_*R\iota^! Rf^!\to Rf^!$ and 
$Rf_!\to Rf_!\iota_*\iota^*$ together with the isomorphism 
$\iota^*\Z/m(n)\cong \Z/m(n)$ induce maps compatible with the pairings
\begin{equation}\label{plko}
\begin{CD} 
H^{i}_c(\X,\Z/m(n))@.\times  H_{i-1}^\et(\X,\Z/m(n))
@>>>  \Z/m\\
@V\iota^* VV @A \iota_* AA 
@|\\
H^{i}_c(Y,\Z/m(n))@. \times H_{i-1}^\et(Y,\Z/m(n))
@>>>  \Z/m.
\end{CD}
\end{equation}

Let $\X$ be connected of dimension $d+1$, smooth and proper over a 
finite field. Then since $Rf^! \Z/m(-n)\cong \Z/m(d+1-n)[2d+2]$,
the pairing above can be rewritten and identified with the cup-product
pairing
$$ H^i_\et(\X,\Z/m(n)) \times H^{2d+3-i}_\et(\X,\Z/m(d+1-n))
\to H^{2d+3}_\et(\X,\Z/m(d+1))\cong  \Z/m. $$
Moreover, we have a trace map
$$H^{2d+4}_\et(\X,\Z(d+1))\stackrel{\partial}{\cong} 
H^{2d+3}_\et(\X,\Q/\Z(d+1))\cong \Q/\Z$$
and compatible cup product pairings for any $m$, 
\begin{equation}
\begin{CD}
H^{i}_\et(\X,\Z(n))@.\times @.H^{2d+4-i}_\et(\X,\Z(d+1-n)) 
@>>>\Q/\Z \\
@VVV @.@A\partial AA @A\partial A\cong A \\
H^{i}_\et(\X,\Z/m(n))@.\times @. H^{2d+3-i}_\et(\X,\Z/m(d+1-n)) @>>> \Q/\Z
\end{CD}\end{equation}

If $X_v$ is connected of dimension $d$, smooth and proper over a 
henselian valuation field of residue characteristic $p$, then the trace
isomorphism $H^{2d}_\et(X_v^\sep,\Q/\ZZ(d))\cong \Q/\ZZ$ together
with the Hochschild-Serre spectral sequence induce isomorphisms
$$H^{2d+3}_\et(X_v,\ZZ(d+1))\stackrel{\partial}{\cong} 
H^{2d+2}_\et(X_v,\Q/\ZZ(d+1))
\cong \Q/\ZZ,$$
and compatible cup product pairings for any $m$ not divisible by $p$,  
\begin{equation}\label{baba}
\begin{CD}
H^{i}_\et(X_v,\ZZ(n))@.\times @.H^{2d+3-i}_\et(X_v,\ZZ(d+1-n)) 
@>>>\Q/\ZZ \\
@VVV @.@A\partial AA @A\partial A\cong A \\
H^{i}_\et(X_v,\Z/m(n))@.\times @. H^{2d+2-i}_\et(X_v,\Z/m(d+1-n)) @>>> \Q/\Z
\end{CD}\end{equation} 
Again the lower pairing is perfect pairing. 

We end this section with a Lemma on topological groups.
\begin{lemma}\label{pdual}
Let $A_i$ be a system of compact groups. Then the natural map  of discrete groups
$$\colim (A_i^*) \to ( \lim A_i)^*$$
is an isomorphism. 
\end{lemma}

The Lemma includes the statement $\oplus (A_i^*) \cong ( \prod A_i)^*$, and applies
in particular to finite groups $A_i$.

\proof 
By Pontrjagin duality, it suffices to prove that the map 
$((\lim A_i)^*)^* \to \lim ((A_i^*)^*)$ obtained by dualizing one
more time is an isomorphism. But since $A_i$ and $\lim A_i$ are compact,
both sides agree with $\lim A_i$.
\endproof

\section{The local situation}
\subsection{Motivic cohomology of the model}
Let $R$ be an excellent henselian discrete
valuation ring with residue field $k$ and field of fractions $K$ (of arbitrary
characteristic). Let $f:\Y\to \Spec R$ be a scheme of finite type, 
and let $m$ be a integer not divisible by $\chr  k$. 
As pointed out by Kahn \cite{kahnhandbook}, the method of \cite{ichmarcII}
and \cite[12.3]{levine} can be
applied to construct a cycle map 
$$cl:  CH_1(\Y,a,\Z/m) \to   H_{a+2}^\et(\Y,\Z/m(1))
:=H^{-a}_\et(\Y, Rf^!\Z/m)$$
which is compatible with localization sequences.
The difference in twists stems from the fact that the left hand side
uses absolute dimension of cycles, and the right hand side the
relative dimension. The map $cl$ induces
a map between the coniveau spectral sequences 
$$ \tilde E^1_{a,b}(\Y,\Z/m)=\oplus_{x\in \Y_{(a)}} 
H^{a-b}_\M(k(x) ,\Z/m(a-1))
\Rightarrow CH_1(\Y,a+b-2,\Z/m)$$
$$ E^1_{a,b}(\Y,\Z/m)=\oplus_{x\in \Y_{(a)}} H^{a-b}(k(x) ,\Z/m(a-1))
\Rightarrow H_{a+b}^\et(\Y,\Z/m(1)),$$
where the $E^1$-terms are motivic cohomology and Galois cohomology, respectively.
By the Rost-Voevodsky theorem, 
$\tilde E^1_{a,b}(\Y,\Z/m)\cong  E^1_{a,b}(\Y,\Z/m)$ for
$b\geq 1$, and $\tilde E^1_{a,b}(\Y,\Z/m)=0$ for $b\leq 0$ for 
cycle dimension reasons. Taking the colimit over powers of $l$, 
the terms  $E^1_{a,b}(\Y,\Q_l/\Z_l)$
vanish for $b<0$ and $l$ a prime different from $\chr k$
by \cite[Lemma 2.6]{saisa}. Hence the difference between
the two theories is measured by the homology of the $E^1_{*,0}$-row, 
$$\cdots\to 
\oplus_{x\in \Y_{(1)}} H^{1}(k(x),\Q_l/\Z_l(0))\to  \oplus_{x\in \Y_{(0)}} 
H^{0}(k(x),\Q_l/\Z_l(-1)). $$
It is shown in \cite[Thm. 2.5.10]{jss} that this complex is isomorphic
up to sign to the complex \cite[Def. 2.1]{saisa} defining Kato homology
$KH_a(\Y,\Q_l/\Z_l)$ so that we obtain a long exact sequence 
\begin{equation}\label{sse} 
\to  KH_{a+3}(\Y,\Q_l/\Z_l) \to  CH_1(\Y,a,\Q_l/\Z_l)\to 
H_{a+2}^\et(\Y,\Q_l/\Z_l) \to KH_{a+2}(\Y,\Q_l/\Z_l)\to .
\end{equation}
Consider the following strengthening of Conjectures 
\cite[Conj. 2.11, 2.12]{saisa} of Saito-Sato:

\begin{conjecture}
Let $\Y$ be a regular scheme, flat and proper over $B$.

1) If $k$ is separably closed, then $KH_a(\Y,\Q_l/\Z_l)=0$ for all $a$.

2) If $k$ is finite, then 
$$KH_a(\Y,\Q_l/\Z_l)=
\begin{cases}
0  &\text{for } a\not=1;\\
(\Q_l/\Z_l)^{I} &\text{for }a=1,
\end{cases} $$
where $I$ is the number of irreducible components of the special fiber 
$\Y\times_Bk$.
\end{conjecture}

\begin{theorem}\label{main1}
The conjecture holds for $a\leq 3$.
\end{theorem}

This is proved in loc.cit. \cite[Thm. 2.13]{saisa} 
for $\Y$ projective such that the reduced special fiber is a strict normal 
crossing scheme. We use
theorems of Gabber and Kerz-Saito to reduce to this case.
Recall that an $l'$-alteration is a proper surjective map, generically
finite of degree prime to $l$.
We will need the following theorem of Gabber 
\cite[Thm. 3 (2)]{gabberthm}\label{gabthm}.

\begin{theorem}
For $B, \Y$ and $l$ as above, there 
exists a finite extension $K'/K$ of degree prime to $l$,
a projective $l'$-alteration $h : \tilde \Y\to  \Y$ over $B'\to B$,
where $B'$ is the normalization of $B$ in $K'$, with
$\tilde \Y$ regular and projective over $B'$.
Moreover, for every geometric point in the closed fiber, 
$\tilde \Y$ is locally for the  etale topology isomorphic to 
$$ B'[t_1,\ldots ,t_n, u_1^{\pm 1},\ldots u_s^{\pm 1}]/
(t_1^{a_1}\cdots t^{a_r}_r u_1^{b_1}\cdots u_s^{b_s}-\pi)$$
at the point $u_i = 1, t_j = 0$, with $1\leq r\leq n$, 
for positive integers $a_1,\ldots, a_r, b_1,\ldots ,b_s$ 
such that $gcd(p,a_1,\ldots , a_r,b_1,\ldots,b_s) = 1$, 
for $p$ the exponential characteristic of $\eta$, and $\pi$
a local uniformizer at $s'$.
\end{theorem}

\begin{remark}
This implies that the reduced special fiber is a normal crossing divisor. 
However, as pointed out in de Jong \cite[\S 2.4]{dejong}, given a normal 
crossing divisor $D\subset \tilde \Y$, there is a projective birational morphism 
$\varphi: \Y'\to \tilde \Y$ 
such that $\varphi^{-1}(D)_{\red}$ us a strict normal crossing divisor.
\end{remark}

\proof[Proof of Theorem \ref{main1}]
The theorem was proved by Saito-Sato 
\cite[Thm.0.8]{saisa} in case that the special fiber
has simple normal crossings, i.e., if its
irreducible components $D_i$ are regular, and if 
the scheme-theoretical intersection $\cap_{i\in J} D_i$ is empty or
regular of codimension $|J|$ in $S$.

In the general case and $a=2,3$, 
we use Gabber's theorem and the remark to find an $l'$-alteration $f: \Y'\to \Y$, 
such that the special fiber 
is a strict normal crossing divisor. By Kerz-Saito \cite[Thm. 4.2, Ex. 4.7]{KS} 
there is a pull-back map 
$f^*: KH_a(\Y,\Q_l/\Z_l)\to KH_a(\Y',\Q_l/\Z_l)$ such that the 
composition $f_*f^*$ with the push-forward is multiplication by the degree of the
alteration, hence the vanishing of $KH_a(\Y',\Q_l/\Z_l)$ implies
the vanishing of $KH_a(\Y,\Q_l/\Z_l)$.

If $a\leq 1$, then by the proper base-change theorem and \cite[(1.9)]{saisa}
we obtain for $d$ the relative dimension of $f$ that 
\begin{multline*}
KH_a(\Y,\Q_l/\Z_l)\cong H^{2-a}_\et(\Y,Rf^!\Q_l/\Z_l)\\
\cong 
H^{2d+2-a}_\et(\Y,\Q_l/\Z_l(d))\cong H^{2d+2-a}_\et(\Y\times_Bk,\Q_l/\Z_l(d)).
\end{multline*}
The latter group vanishes for $a\leq 1$ and $k$ separably closed,
$a \leq 0$ and $k$ finite, and is isomorphic to $(\Q_l/\Z_l)^I$ for
$a=1$ and $k$ finite \cite[Lemma 2.15]{saisa}.
\endproof 

From the sequence \eqref{sse} we immediately obtain:

\begin{corollary} If $k $ is finite, we have a surjection 
$CH_1(\Y,1,\Q/\ZZ)\twoheadrightarrow H_{3}^\et(\Y,\Q/\ZZ(1)),$
isomorphisms
$CH_1(\Y)\otimes\Q/\ZZ\cong H_2^\et(\Y,\Q/\ZZ(1))$ and 
$H_{1}^\et(\Y,\Q/\ZZ(1)) \cong (\Q/\ZZ)^I$, and 
$H_i^\et(\Y,\Q/\ZZ(1))$ vanishes for $i\leq 0$.
\end{corollary}

We give an application to integral  etale motivic cohomology, which can 
be viewed as a generalization of the vanishing of Brauer groups of
relative curves over $B$ \cite[\S 3]{brauerIII}. We return to 
cohomological notation $H^i_\et(\Y,\ZZ(d))= H_{2d+2-i}^\et(\Y,\ZZ(1))$
for regular $\Y$ of dimension $d+1$.

\begin{corollary}\label{descc1}
Let $\Y$ be a regular scheme of pure dimension $d+1$, flat and proper over $B$
with finite residue field of characteristic $p$.
Then we have a surjection
$$CH_1(\Y)_{\Z'} \twoheadrightarrow H^{2d}_\et(\Y,\ZZ(d)),$$ 
and isomorphisms
$H^{2d+1}_\et(\Y,\ZZ(d)) \cong 0$ and $H^{2d+2}_\et(\Y,\ZZ(d))\cong (\Q/\ZZ)^I$.
\end{corollary}

\proof
In degrees $\leq 2d+1$, 
this follows by combining the Corollary, the isomorphism
$H^i_\et(\Y,\Q(d))\cong CH_1(\Y,2d-i)_\Q$ (which vanishes for $i>2d$),  
and the coefficient sequence
$$\begin{CD} 
CH_1(\Y,1)_\Q@>>> CH_1(\Y,1,\Q/\ZZ) @>>> CH_1(\Y)_{\Z'} @>>> CH_1(\Y)_\Q \\
@|@VsurjVV @VVV @|  \\
H^{2d-1}_\et(\Y,\Q(d))@>>> H^{2d-1}_\et(\Y,\Q/\ZZ(d)) @>>> 
H^{2d}_\et(\Y,\ZZ(d)) @>>> H^{2d}_\et(\Y,\Q(d))
\end{CD}$$ $$ \begin{CD}
@>>> CH_1(\Y)\otimes \Q/\ZZ @>>> 0\\
@. @| @VVV  \\
@>>> H^{2d}_\et(\Y,\Q/\ZZ(d)) @>>> H^{2d+1}_\et(\Y,\ZZ(d)) @>>> 0.
\end{CD}$$
Finally, 
$H^{2d+2}_\et(\Y,\ZZ(1)) \cong H^{2d+1}_\et(\Y,\Q/\ZZ(1))\cong (\Q/\ZZ)^I $. 
\endproof

\subsection{Motivic Cohomology of the generic fiber}
Let $K$ be a global field of exponential
characteristic $p$ and $X$ be a connected scheme of dimension $d$, 
smooth and projective over $K$. For a valuation
$v$ of $K$, we let $K_v$ be the henselization of $K$ at $v$, 
$G_v$ its Galois group, and $X_v=X\times_KK_v$. 

\begin{proposition}\label{highdeg}
We have
$$H^{2d+2}_\et(X_v,\ZZ(d))\cong H^{2}(K_v,\ZZ)$$
and 
$$H^{2d+2}_\et(X,\ZZ(d))\cong H^2(K,\ZZ). $$
The map $l_d^{2d+2}$ is injective. 
\end{proposition}

\proof 
Since $X_v$ has dimension $d$, Prop. \ref{rojtman} and 
the Hochschild-Serre spectral sequence gives
$$H^{2d+2}_\et(X_v ,\ZZ(d))\cong H^2(K_v,H^{2d}_\et(X_v^\sep,\ZZ(d)))
\cong H^2(K_v, CH_0(X^\sep_v)_{\Z'}).$$ 
By Proposition \ref{rojtman}, we have 
$H^2(K_v, CH_0(X^\sep_v)_{\Z'}^0)\cong H^2(K_v,\Alb_{\Z'})$,
and the latter group vanishes by comparing to the completion of $K_v$
and using \cite[Rem. I 3.6, 3.10]{adt}. 
Hence  the degree map induces an isomorphism
$$H^2(K,H^{2d}_\et(X^\sep_v,\ZZ(d)))\cong H^2(K_v,\ZZ).$$

In the global case, the same argument works because $H^2(K,\Alb_{\Z'})$ still 
vanishes by \cite[I Cor. 6.24]{adt}.
The injectivity of $l_d^{2d+2}:H^1(K,\Q/\ZZ)\to \prod H^1(K_v,\Q/\ZZ)$ 
follows from Chebotarev's density theorem.
\endproof

The following is a refinement of a theorem of Saito-Sato, 
see \cite[Thm. 3.25]{ct}:

\begin{theorem}\label{ss2}
Let $X$ be a smooth and proper over the fraction
field of an excellent henselian discrete valuation ring with finite or 
separably closed residue field. Then 
the group of zero-cycles of degree $0$ on $X$ is isomorphic to the direct sum
of a finite group and a group divisible by all integers prime to the 
residue characteristic.
\end{theorem}

We use this to prove the following 

\begin{theorem} \label{localthm}
Let $X$ be smooth and proper of dimension $d$ 
over a global field of characteristic $p$.
Then the cokernel $B_v$ of
$H^{2d}_\et(X_v,\ZZ(d))\to H^{2d}_\et(X^\sep_v,\ZZ(d))^{G_v}$
is finite and vanishes for almost all $v$.
\end{theorem}

\proof
We show that the cokernel is finite, and vanishes for $X_v$ with good reduction. Consider the map of exact sequences
$$\begin{CD}
0@>>> H^{2d}_\et(X_v,\ZZ(d))^0@>>> H^{2d}_\et(X_v,\ZZ(d))@>\deg>> \ZZ 
@>>> \ZZ/\delta_v\ZZ@>>> 0\\
@.@VVV @VVV@| @VVV \\
0@>>> (H^{2d}_\et(X^\sep_v,\ZZ(d))^0)^G @>>> H^{2d}_\et(X^\sep_v,\ZZ(d))^G 
@>\deg>> \ZZ^G @>>> \ZZ/\delta'_v\ZZ@>>> 0,
\end{CD}$$
where the map $\deg$ is induced by the proper push-forward along the structure
morphism \cite[Cor. 3.2]{ichduality}. 
The invariants $\delta_v$ and $\delta'_v$ are 
analogs of the (prime to $p$-part of the)
index and period of $X_v$. The cokernels of
the two left vertical maps differ by a finite group, and 
they agree if $X_v$ has good reduction because in this case 
a zero-cycle of degree $1$ in the special fiber 
can be lifted to $X_v$ by the henselian property,
hence $\delta_v=\delta'_v=1$.
Now consider the following diagram with exact rows 
(we omit the coefficients $\ZZ(d)$):
$$\begin{CD}
0@>>> \Tor\; CH_0(X_v)^0_{\Z'} @>>>  CH_0(X_v)^0_{\Z'} @>\tau_1>>  
CH_0(X_v)^0_\Q \\
@. @VVV @VVV @| \\
0@>>> \Tor H^{2d}_\et(X_v)^0 @>>> H^{2d}_\et(X_v)^0 @>\tau_2>> 
H^{2d}_\et(X_v)^0_\Q \\ 
@. @VVV @VVV @| \\ 
0@>>> (\Tor H^{2d}_\et(X^\sep_v)^0)^{G_v} @>>> 
(H^{2d}_\et(X^\sep_v)^0)^{G_v}  @>\tau_3>> 
(H^{2d}_\et(X^\sep_v)^0\otimes \Q)^{G_v} .
\end{CD}$$
The lower right equality sign follows by a trace argument.
By Theorem \ref{ss2}, 
the cokernel $CH_0(X_v)^0\otimes\Q/\ZZ$ of $\tau_1$ vanishes. 
This implies that $\tau_2$ and $\tau_3$ are surjective,  
and the cokernels of the two lower left vertical maps are isomorphic.
By Proposition \ref{rojtman}  we get
$$(\Tor H^{2d}_\et(X^\sep_v,\ZZ(d))^0)^{G_v}
\cong (\Tor \Alb_{X_v}(K^\sep_v)_{\Z'} )^{G_v}
\cong \Tor \Alb_{X_v}(K_v)_{\Z'}.$$
This group, hence $B_v$, is finite, because an abelian variety over a 
local field
has a subgroup of finite index which is divisible by all primes different
from $p$ \cite[I Rem. 3.6]{adt}.

\begin{proposition}
Let $X_v$ be a smooth and proper scheme with good reduction over the 
henselian valuation field $K_v$. Then the map 
$$ H^i_\et(X_v,\Q/\ZZ(n))\stackrel{\sigma }{\to }
H^i_\et(X_v^\sep,\Q/\ZZ(n))^{G_v}$$
is surjective. 
In particular, the map 
$$\Tor  H^i_\et(X_v,\ZZ(n))\stackrel{\sigma }{\to }
(\Tor H^i_\et(X_v^\sep,\ZZ(n)))^{G_v}$$
is surjective for $i\not=2n+1$. 
\end{proposition}

This implies the vanishing of $B_v$ for $v$ with good reduction
because  the torsion of $H^{2d}_\et(X_v^\sep,\ZZ(d))$ and
of $H^{2d}_\et(X_v^\sep,\ZZ(d))^0$ agree.

\proof
Let $Y$ be the special fiber of a smooth and proper model 
$\Y$ over the valuation ring $R$ of $K_v$ with strict 
henselization $R^{nr}$
and residue field $k$. Consider the specialization diagram
$$\begin{CD}
H^i_\et(Y,\Z/m(n)) @<\sim << H^i_\et(\Y,\Z/m(n)) 
@>>> H^i_\et(X_v,\Z/m(n))\\
@VVV @VVV @VVV \\
H^i_\et(Y^\sep,\Z/m(n))^{\Gal(k)} @<\sim << 
H^i_\et( \Y\times_RR^{nr},\Z/m(n))^{\Gal(k)}  
@>\sim >> H^i_\et(X_v^\sep,\Z/m(n))^{G_v}.
\end{CD}$$
The left horizontal maps are isomorphisms by the proper base change theorem, and
the lower right map is an isomorphism by the smooth base change theorem.
The first statement follows because $k$ has cohomological dimension one,
hence the left vertical map is surjective.
The second statement follows because 
$H^{i-1}_\et(X_v^\sep,\Q/\ZZ(n))\cong \Tor H^i_\et(X_v^\sep,\ZZ(n))$
for $i\not=2n+1$ by Proposition \ref{strstr}. 
\endproof

\section{Proof of \ref{order}, \ref{maines}, \ref{secondes}}
Let $C$ is a smooth and proper 
curve over a finite field $\F$ of characteristic $p$, with generic point 
$\Spec K\to C$. Let $\X\to C$ be a
flat, projective map, with regular $\X$ of dimension $d+1$. 
We assume that the generic fiber $X=\X\times_CK$ is smooth 
and geometrically connected over $K$.
For a closed point $v$ of $C$, we let $\O_v$ be the henselization of 
$C$ at $v$, $K_v$ its quotient field, 
$k_v$ be the residue field, and $\X_v=\X\otimes_C\O_v$,  
with generic fiber $X_v=X\times_KK_v$ and closed fiber $Y_v=\X\times_C k_v$.
Let $G$ be the Galois group of $K$ and $G_v$ the Galois groups of $K_v$.

\subsubsection*{Proof of Proposition \ref{secondes}(2)}
Taking the colimit over increasing union of fibers of the localization sequence
in etale cohomology, we obtain a map of long exact sequences
\begin{equation*}
\begin{CD}@>>> H^{i}_\et(\X,\Z(n))@>>>  
H^{i}_\et(X,\Z(n)) @> \partial>>  
\bigoplus_v H^{i+1}_{Y_v}(\X,\Z(n)) @>>> \\
@. @VVV  @Vl^i_nVV @VVV \\
@>>> \prod_v H^{i}_\et(\X_v,\Z(n))@>>>  
\prod_v H^{i}_\et(X_v,\Z(n)) @> (\partial_v)>>  
\prod_v H^{i+1}_{Y_v}(\X_v,\Z(n))@>>>
\end{CD}\end{equation*}
We claim that identifying the terms $H^{i}_\et(\X_v,\Z(n))$ with 
Corollary \ref{descc1} this gives rise to the following 
commutative diagram. 
\begin{equation}\label{mmdd}\begin{CD}
H^{2d+1}_\et(\X,\ZZ(d))@>>> \phantom{\bigoplus_v} H^{2d+1}_\et(X,\ZZ(d)) @>\partial>> 
\bigoplus_v H^{2d+2}_{Y_v}(\X,\ZZ(d))\\
@VVV @Vl^{2d+1}_d VV  @| \\
0@>>> \bigoplus_v H^{2d+1}_\et(X_v,\ZZ(d)) @>(\partial_v)>> 
\bigoplus_v H^{2d+2}_{Y_v}(\X_v,\ZZ(d))
\end{CD}\end{equation}
$$\begin{CD}@>>> H^{2d+1}_\et(\X,\Q/\ZZ(d))@>>>  H^{2d+1}_\et(X,\Q/\ZZ(d)) \\
@. @V\delta VV @Vl^{2d+2}_dVV \\
@>(\tau_v)>> \prod_v \Q/\ZZ^{I_v} @>>>\prod_v H^{2d+1}_\et(X_v,\Q/\ZZ(d)).
\end{CD}$$
The shift of degrees appears 
because we switch from $\Z'$ to $\Q/\ZZ$-coefficients
in the last four terms. The map 
$l_d^{2d+1}$ has image in the direct sum because $\partial$ has.
The lower sequence is exact because for almost all $v$, 
$\tau_v$ is the zero map. Indeed, if $\X_v$ has good reduction, 
then by the proper base change theorem 
$H^{2d+1}_\et(\X_v,\Q/\ZZ(d))\cong H^{2d+1}_\et(Y_v,\Q/\ZZ(d))
\cong H^{1}_\et(k_v,\Q/\ZZ)$, so that using Proposition \ref{highdeg}
the map $H^{2d+1}_\et(\X_v,\Q/\ZZ(d))\to  H^{2d+1}_\et(X_v,\Q/\ZZ(d))$
can be identified with the map 
$H^{1}_\et(k_v,\Q/\ZZ)\to  H^{1}_\et(K_v,\Q/\ZZ)$, which is the injection
dual to the surjection $\Gal(K_v)\to \Gal(k_v)$.

From the diagram and Proposition \ref{highdeg} we see that  
$\ker l^{2d+1}_d=\im H^{2d+1}_\et(\X,\ZZ(d))$, and  
$\coker  l^{2d+1}_d\cong \ker \delta$,
where $\delta$ is the colimit of the maps  
$$\delta_m: H^{2d+1}_\et(\X,\Z/m(d))\stackrel{\iota^*}{\to} 
\prod_v H^{2d+1}_\et(Y_v,\Z/m(d))
\stackrel{tr}{\to } \prod_v (\Z/m)^{I_v}$$
followed by the natural inclusion 
$\colim_{p\not|m} \prod_v (\Z/m)^{I_v}\subseteq \prod_v (\Q/\ZZ)^{I_v}$.

\begin{proposition}\label{cokerl}
Up to $p$-groups, we have an exact sequence
$$0\to (T\Br(\X))^* \to \ker  \delta \to \Pic(X)^* \to 0$$
\end{proposition}

\proof 
By Poincar\'e duality and Lemma \ref{pdual} we have 
$$
 H^{2d+1}_\et(\X,\Q/\ZZ(d))
\cong \colim_{p\not|m} H^2_\et(\X,\Z/m(1))^* 
\cong H^2_\et(\X,\hat\ZZ(1))^*,$$
where here and in the following $H^i_\et(\X,\hat\ZZ(n))$
denotes $\lim_{p\not| m} H^i_\et(\X, \Z/m(n))$.
If $\iota :Y_v \to \X$ is the inclusion of a closed fiber, 
then we obtain the commutative diagram \eqref{plko}
$$\begin{CD} 
H^{2d+1}_\et(\X,\Z/m(d))@.\times  H^2_\et(\X,\Z/m(1))@>>>
\Z/m\\
@V\iota^* VV @A \iota_* AA 
@|\\
H^{2d+1}_\et(Y_v,\Z/m(d))@.\times H_{2d}^\et(Y_v,\Z/m(d))@>>> 
\Z/m.
\end{CD}$$
Thus $\delta_m$ 
is dual to the right vertical map in
$$\begin{CD}
\bigoplus_v \CH_{d}(Y_v)/m @= \bigoplus_v H_{2d}^\et(Y_v,\Z/m(d))\\
@V\iota_*^mVV @Vd_mVV \\
\Pic(\X)/m @>\subset >>  H^2_\et(\X,\Z/m(1)).
\end{CD}$$
Since $\coker \iota_*^m=\Pic(X)/m$, the coefficient sequence gives 
a short exact sequence of cokernels 
$$ 0\to \Pic(X)/m \to \coker d_m \to {}_m\Br\X\to 0$$
and its dual 
$$ 0\to ({}_m\Br\X)^*  \to \ker \delta_m \to (\Pic(X)/m)^* \to 0.$$
In view of the finiteness of the groups involved, the proposition follows by taking
colimits using Lemma \ref{pdual}.
\endproof

\begin{lemma}\label{whatmap}
The composition 
$a: \coker l_d^{2d+1}\stackrel{\sim}{\to} \ker \delta\to\Pic(X)^*$
is induced by the composition
$a': H^{2d+1}_\et(X_v,\ZZ(d))\to \Pic(X_v)^* \to \Pic(X)^*$,
where the first map is the duality \eqref{baba}.
\end{lemma}

\proof
Consider the following diagram, where the upper horizontal maps
are induced by \eqref{mmdd}, and the lower horizontal maps
are functoriality.
$$\xymatrix{
H^{2d+1}_\et(X_v,\ZZ(d)) \ar[r]^{\partial}\ar[d]&
\oplus H^{2d+2}_{Y_v}(\X,\ZZ(d))\ar[r]\ar[d]&
H^{2d+2}_\et(\X,\ZZ(d))\ar[d]\\
\Pic(X_v)^*\ar[r]\ar@/_/[dr] &
\Pic(\X_v)^* \ar[r]&
\Pic(\X)^*\\
& \Pic(X)^*\ar@/_/@{^{(}->}[ur]}
$$
The upper vertical maps are the cup products over local fields, over
local rings, and over finite fields.
Commutativity amounts to the compatibility of these pairings.
The lower triangle is commutative by functoriality.
The upper horizontal composition induces the identification of 
$\coker l_d^{2d+1}$ with $\ker \delta$ and the right vertical map
induces the map $\ker \delta\to \Pic(X)^*$ of the Proposition. 
Commutativity of the diagram shows that $a$ is induced by the lower left 
composition.
\endproof

In particular, we obtain Proposition \ref{secondes} because 
$\coker l_d^{2d+1}\cong \ker\delta\cong \Pic (X)^*$ if $\Br(\X)$ is finite.

\subsubsection*{First diagram}
We consider the long exact sequence of Galois cohomology groups
associated to the degree map 
$$0\to CH^d(X^\sep)^0_{\Z'}\to 
CH^d(X^\sep)_{\Z'}\stackrel{\deg}{\longrightarrow} \ZZ\to 0.$$
Since higher Galois cohomology of a uniquely divisible groups vanish, 
Proposition \ref{rojtman} gives an isomorphism 
$H^1(K, CH^d(X^\sep )^0_{\Z'})\cong H^1(K,\Alb_X)_{\Z'}$. 
Let 
$D$ and $D_v$ be the cokernels of the global and local degree maps 
$$
CH^d(X^\sep)^G_{\Z'}\to \ZZ, \qquad
CH^d(X_v^\sep)^{G_v}_{\Z'}\to \ZZ$$
and recall that their orders were denoted by $\delta'$ and $\delta'_v$. 
Then we obtain the following commutative diagram with exact middle rows.
The map $\beta^1$, hence $\tau$, has image in the direct sum
by \cite[I Lemma 6.3]{adt}.
\begin{equation}\label{dia2}
\begin{CD}
@.K_2@>>> \Sh(\Alb_X) @>>> \Phi\\
@.@VVV @VVV @VVV\\
0@>>> D@>>> H^1(K,\Alb_X )_{\Z'}@>>> 
 H^1(K,H^{2d}_\et(X^\sep,\ZZ(d)))@>>> 0 \\
@.@V\theta VV @V\beta^1 VV @V\tau VV \\
0@>>> \bigoplus_v D_v@>>>
\bigoplus_v H^1(K_v,\Alb_X)_{\Z'}@>u>> 
\bigoplus_v H^1(K_v,H^{2d}_\et(X^\sep,\ZZ(d)))@>>> 0\\
@.@VVV @VVV @V\rho VV \\
@.C_2@>>> (T\Sel(\Pic^0_X))^*@>>> \Psi.
\end{CD}\end{equation}
The groups $K_2$ and $C_2$ are the kernel and cokernel
of $\theta$, and $\Phi$ and $\Psi$ are the kernel and cokernel of $\tau$, 
respectively. 
The outer columns are exact by definition, 
and the middle column is exact because for any abelian variety $A$ with dual 
$A^t$ there is a short exact sequence \cite{tan}
$$0\to  \Sh(A)\to  H^1(K,A )\stackrel{\beta}{\longrightarrow}
\bigoplus_v H^1(K_v,A) \to  (T\Sel(A^t))^*\to  0 $$
with
$T\Sel(A^t)=\lim_{p\not|m} \ker H^1(K,{}_mA^t)\to \prod_v H^1(K_v,A^t)$
the inverse limit of Selmer groups.
The middle two rows induce the horizontal 
maps in the upper and lower row. 
An easy diagram chase in \eqref{dia2} gives the exact sequence
\begin{equation}\label{iii}
 0\to K_2\to \Sh(\Alb_X)\to \Phi \to C_2\to (T\Sel(\Pic^0_X))^* 
\to \Psi\to 0
\end{equation}

\begin{proposition}\label{erste}
1) The maps $\Sh(\Alb_X)\to \Phi$
and $(T\Sel(\Pic^0))^* \to \Psi$ have finite kernels and cokernels. 

2) If $\Sh(\Alb_X)$ is finite, then $\Phi$ is finite, 
$u$ induces a surjection $z: H^0(K,\Pic^0_X)^*\to \Psi$, 
and we have an equality 
\begin{equation}\label{for1}
\delta' \cdot |\Phi|\cdot |\ker z| 
=|\Sh(\Alb_X)|\cdot \prod \delta'_v.
\end{equation}
\end{proposition}

\proof 
1) Follows because $D$ and $\prod_vD_v$,  
hence $K_2$ and $C_2$ are finite.

2) If $\Sh(A)$ is finite, then $T\Sh(A)=0$ and hence
$A(K)^\wedge \cong T\Sel(A)$. Thus $(T\Sel(\Pic^0_X))^* \cong 
H^0(K,\Pic^0_X)^*$ and we obtain the map $z$. 
Now the equality follows by taking alternating orders in \eqref{iii}
because $|K_2|/|C_2|=
\delta'/\prod_v\delta_v'$.
\endproof

\subsubsection*{The second diagram}
Let $B$ and $B_v$ be the cokernels of 
$$H^{2d}_\et(X,\ZZ(d))\stackrel{\rho_d}{\longrightarrow} H^{2d}_\et(X^\sep,\ZZ(d))^G,\quad
H^{2d}_\et(X_v,\ZZ(d))\stackrel{}{\longrightarrow}  H^{2d}_\et(X_v^\sep,\ZZ(d))^{G_v},$$
respectively, and recall that we denoted their orders by $\beta$ and $\beta_v$.  
Then $\prod_v\beta_v$ is finite by Theorem \ref{localthm}.
From the Hochschild-Serre spectral sequence we obtain the sequence
\begin{multline*}
0\to B\to  H^2(K,H^{2d-1}_\et(X^\sep,\ZZ(d)))\\ \to 
H^{2d+1}_\et(X,\ZZ(d)) \to
H^1(K,H^{2d}_\et(X^\sep,\ZZ(d)))\to 0. 
\end{multline*}
By Proposition \ref{strstr},  
the injection
$$ H^{2d-2}_\et(X^\sep,\Q/\ZZ(d)) \cong 
\Tor H^{2d-1}_\et(X^\sep,\ZZ(d))) \hookrightarrow
H^{2d-1}_\et(X^\sep,\ZZ(d))) $$
has uniquely divisible cokernel, hence we obtain an isomorphism
$$H^2(K,H^{2d-2}_\et(X^\sep,\Q/\ZZ(d)))\cong 
H^2(K,H^{2d-1}_\et(X^\sep,\ZZ(d))).$$

Comparing this with the same exact sequence for the local situation,
we obtain the exact middle two rows of the following diagram.

{
\footnotesize
\begin{equation}\label{eins+}
\begin{CD} 
K_1 @>>> \ker \xi_d^{2d-2} @>>> \ker l_d^{2d+1}@>>> \Phi \\
@VVV @VVV @VVV @VVV \\
B @>>> H^2(K,H^{2d-2}(X^\sep,\Q/ \ZZ(d))) @>>> H^{2d+1}_\et(X,\ZZ(d))@>>> 
H^1(K,H^{2d}_\et(X^\sep,\ZZ(d)))\\
@V\theta' VV @V\xi_d^{2d-2} VV @Vl^{2d+1}_dVV @V\tau VV \\
\bigoplus B_v@>>> \bigoplus H^2(K_v,H^{2d-2}(X^\sep,\Q/\ZZ(d)))@>>>\bigoplus  
H^{2d+1}_\et(X_v,\ZZ(d))@>>> 
\bigoplus H^1(K_v,H^{2d}_\et(X_v^\sep,\ZZ(d)))\\
@VVV @V\xi'VV @VVV@V\rho VV \\
C_1@>>> (H^2_\et(X^\sep,\hat\Z'(1))^G)^* @>w>>  \coker l_d^{2d+1} @>>> \Psi
\end{CD}
\end{equation}}
The maps $\xi_d^{2d-2}$ and $\tau$ have image in the direct sum by 
\cite[I Lemma 4.8, 6.3]{adt}, and we have seen above that
$l_d^{2d+1}$ has image in the direct sum. The groups $K_1$ and $C_1$
are the kernel and cokernel of the canonical map $\theta'$, and 
the map $\xi'$ is the dual of the injection
$H^2_\et(X^\sep,\Z/m(1))^G\to \prod_v H^2_\et(X^\sep,\Z/m(1))^{G_v}$.
All columns are short exact by definition except  the left middle column,
which is exact by Tate-Poitou duality. 
Hence the maps in the middle two rows induce the maps in the upper
and lower row.

\subsubsection*{Proof of Theorem \ref{maines}.}
In view of the finiteness of 
$\prod_v \beta_v$ and Proposition \ref{erste}, diagram
\eqref{eins+} gives a sequence defined and exact up to 
finite groups and $p$-groups
\begin{multline*}
0\to  \coker \rho_d \to  \ker \xi_d^{2d-2} \to  \ker l_d^{2d+1}
\to  \Sh(\Alb_X)\to \\
(H^2_\et(X^\sep,\hat\Z'(1))^G)^* \stackrel{w}{\to}  
\coker l_d^{2d+1} \to (T\Sel(\Pic^0_X))^*\to 0.
\end{multline*}

\begin{lemma}
The composition of $w$ with the map $a$ of Lemma \ref{whatmap}
is the dual of the cycle map 
$\Pic(X)^\wedge\to H^2_\et(X^\sep,\hat\Z'(1))^G$.
\end{lemma}

\proof
Consider the following commutative diagram. The left middle square
is commutative because duality of Galois cohomology of a local field 
is compatible with duality of etale cohomology over local fields,
and the commutativity of the part with $a'$ is \eqref{baba}.
$$\xymatrix{
H^2(K_v,H^{2d-1}_\et(X^\sep,\ZZ(d))) \ar[r] &
H^{2d+1}_\et(X_v,\ZZ(d))\ar@/^3pc/[rdd]^{a'}\\
H^2(K_v,H^{2d-2}_\et(X^\sep,\Q/\Z'(d))) \ar[u]^\cong \ar[r]\ar[d]_\cong&
H^{2d}_\et(X_v,\Q/\Z'(d))\ar[d]_\cong \ar@{>>}[u]\\
(H^2_\et(X_v^\sep,\hat\Z(1))^{G_v})^*\ar[d]\ar[r] & 
H^2_\et(X_v,\hat\Z(1))^* \ar[d]\ar[r]&
\Pic(X_v)^*\ar[d]\\
(H^2_\et(X^\sep,\hat\Z(1))^{G})^* \ar[r]&
H^2_\et(X,\hat\Z(1))^*\ar[r] &
\Pic(X)^*}$$
By Lemma \ref{whatmap}, the map $a\circ w$ is induced by lifting an 
element along the left vertical map and taking the upper right composition.
Commutativity implies that this is the map in the lower row, i.e., the cycle map.
\endproof

Since $\Pic^0(X)^\wedge$ is the kernel of the cycle map, 
we obtain a commutative diagram with exact rows and surjective map $c_1^*$.
\begin{equation}\begin{CD}\label{goesto}
@>>>(H^2_\et(X^\sep,\hat\Z'(1))^G)^* @>w>>
\coker l_d^{2d+1} @>>> \Psi @>>> 0\\
@.@Vc_1^*VV @VaVV@VbVV\\
0@>>> \NS(X)^* @>>>\Pic(X)^*@>>> \Pic^0(X)^*@>>> 0
\end{CD}\end{equation}
Taking the kernel of the vertical maps, 
we obtain the sequence of Theorem \ref{maines} (still up to finite groups
and $p$-groups) because $\ker a\cong (T\Br(\X))^*$ by Proposition \ref{cokerl}
and $\ker b\cong (T\Sh(\Pic^0_X))^*$ by Proposition \eqref{erste}(1) and the definition
of the Selmer group (it is not necessary to check that $b$ is the canonical map,
as $b$ is surjective so that $\ker b$ and $(T\Sh(\Pic^0_X))^*$
have the same corank and are abstractly isomorphic).

\subsubsection*{Proof of Theorem \ref{order}.}
If $\Br(\X)$, hence $\Sh(\Alb_X)$,
$\Sh(\Pic^0_X)$ and $\coker c$ are finite, then 
the maps $a$ and  $b$ in \eqref{goesto} are isomorphisms. 
This implies that $\ker w=\ker c_1^*=\ker c^*$, hence 
counting alternating orders of the first and fourth row of \eqref{eins+}, 
we get 
\begin{equation}\label{for2}
|\ker \xi_d^{2d-2}/K_1| \cdot |\Phi| \cdot |C_1|=
|\ker l_d^{2d+1}| \cdot |\coker c| .
\end{equation}
We now compare \eqref{for1} and \eqref{for2}.
We have the identities $|C_1|= \prod\beta_v/|B/K_1| $, and by
the following Lemma 
$|\ker z|=|\ker bz|$ is equal to the order $\alpha$ of the cokernel
the canonical map $\Pic^0(X)\to H^0(K, \Pic^0_X)$, so that
$$|\ker l_d^{2d+1}|\cdot |\coker c|\cdot |B/K_1|\cdot \alpha  \delta'=
|\Sh(\Alb_X)|\cdot |\ker \xi_d^{2d-2}/K_1|\prod_v \beta_v\delta'_v.$$
If $B=\coker \rho_d$ is finite, we can remove the common factor $|K_1|$
and obtain the formula of Theorem \ref{order}. It remains to show:

\begin{lemma} 
If $\Sh(\Pic^0_X)$ is finite, then the composition
$H^0(K,\Pic^0_X)^*\stackrel{z}{\to} \Psi\stackrel{b}{\to} \Pic^0(X)^*$ 
of the maps from Proposition \ref{erste} and diagram \eqref{goesto} 
is the canonical map.
\end{lemma}

\proof
Consider the filtration induced by the Hochschild-Serre spectral sequence. 
By functoriality, we have a duality between 
$$H^2_\et(X_v,\hat\ZZ(1))/F^2=\coker\left( \hat\ZZ \cong 
H^2(K_v,\hat\ZZ(1))\to H^2_\et(X_v,\hat\ZZ(1))\right)$$
and 
$$H^{2d}_\et(X_v,\Q/\ZZ)^0=
\ker\left( H^{2d}_\et(X_v,\Q/\ZZ)\to 
H^{2d}_\et(X_v^\sep,\Q/\ZZ(d))^{G_v}\cong \Q/\ZZ \right).$$
The latter still surjects onto 
$H^{2d+1}_\et(X_v,\ZZ(d))$ 
because $H^{2d}_\et(X_v,\Z(d))\otimes\Q/\ZZ$ surjects onto
the copy of $\Q/\ZZ$ and maps to zero in $H^{2d+1}_\et(X_v,\ZZ(d))$.
We obtain a diagram
\begin{equation}\label{xrxr}
\xymatrix{
H^{2d+1}_\et(X_v,\ZZ(d))\ar@{>>}[r]\ar@/_7pc/[ddd]_{a'}&
H^1(K_v,H^{2d}_\et(X^\sep,\ZZ(d)))\ar@{>>}[r]^-{\rho}&
\Psi\\
H^{2d}_\et(X_v,\Q/\ZZ(d))^0\ar@{>>}[u]\ar@{>>}[r]\ar[d]_\cong&
H^1(K_v,\Alb_{X})\ar@{>>}[u]^u\ar[d]_\cong\\
(H^2_\et(X_v,\hat\ZZ(1))/F^2)^* \ar@{>>}[r]\ar[d]&
H^0(K_v,\Pic^0_{X})^* \ar[d]\ar[r]&
H^0(K,\Pic^0_{X})^*\ar[uu]^z\ar[d]\\
\Pic(X_v)^*\ar[r]\ar@/_/[rd] &
\Pic^0(X_v)^*\ar[r] &
\Pic^0(X)^*\\
& \Pic(X)^* \ar@/_/[ru]
}\end{equation}
The map $u$ is the map of diagram \eqref{dia2}, and 
the upper right square defines the map $z$.  
The map $b$ is induced by lifting elements of $\Psi$ 
to $H^{2d+1}_\et(X_v,\ZZ(d))$
and then following the curved arrows.  
The middle left square is the following composition: 
\begin{multline*}
\begin{CD}
H^{2d}_\et(X_v,\Q/\ZZ(d))^0@>>> 
H^1(K_v,H^{2d-1}_\et(X_v^\sep,\Q/\ZZ(d)))@>\sim>>\\
@V\cong VV @V\cong VV  \\
(H^2_\et(X_v,\hat\ZZ(1))/F^2) ^*@>>> 
H^1(K_v,H^1_\et(X_v^\sep,\hat\ZZ(1)))^*@>\sim>> \\ 
\end{CD}\\
\begin{CD}
@>\sim >>H^1(K_v,\Tor  \Alb_{X})
@>\sim>> H^1(K_v,\Alb_{X})\\
@. @V\cong VV @V\cong VV \\
@<\sim << H^1(K_v, T\Pic^0_{X})^*@>\sim>> H^0(K_v,\Pic^0_{X})^*.
\end{CD}
\end{multline*}
The left square commutes by compatibility of duality for schemes over 
local fields with Galois cohomology of the field, and 
the middle square is the compatibility 
of the $e_m$-pairing with the cup-product pairing:
$$\begin{CD} 
H^{2d-1}_\et(X_v^\sep,\Q/\ZZ(d))@>\sim >> \Tor \Alb(X_v^\sep)\\
@V\cong VV @V\cong Ve_mV \\
H^1_\et(X_v^\sep,\hat \ZZ(1))^*@<\sim<< (T \Pic^0(X_v^\sep))^*.  
\end{CD}$$
It is easy to see that all other squares in \eqref {xrxr} commute,
and this implies the claim of the Lemma.  
\endproof


If $\X$ is a surface, i.e. $d=1$, 
then $\ker l_d^{2d+1}=\Br(\X)$, $\ker \xi_1^0$ vanishes, and 
$\delta'$ and $\delta'_v$ are the periods of $X$ and $X_v$. 
By Lichtenbaum \cite{lichtenbaumcurve}, 
$\beta_v$ is equal to the index $\delta_v$ of $X_v$, whereas 
$c$ is the (completed) degree map whose cokernel has order equal 
to the index $\delta$. 
Hence the formula reduces to the formula in \cite{ichsurface}
in view of $\delta\alpha=\beta\delta'$, which can be seen from the
following diagram:
$$ \begin{CD}
0@>>> \Pic^0(X) @>>> \Pic X@>>> \delta \Z@>>> 0\\
@. @VVV @VVV @VVV \\
0@>>> H^0(K,\Pic^0_{X}) @>>> H^0(K,\Pic_{X})@>>> \delta' \Z@>>> 0
\end{CD}$$

\section{Other consequences}
We recall the following well-known results: 
\begin{enumerate}
\itemsep0em
\item Tate's conjecture on the surjectivity of the cycle map for $\X$ 
in codimension $n$ is equivalent to the finiteness of
$H^{2n+1}_\et(\X,\Z(n))$, or the finiteness of its $l$-primary part for any $l$
\cite[Prop. 3.2]{ichtata}. In particular, Tate's conjecture for divisors on $\X$ holds
if and only if $\Br(\X)$ is finite.
\item The finiteness of the Brauer group and of the Tate-Shafarevich group is implied
by the finiteness of its $l$-primary part for any prime $l$, see 
\cite[Remark 8.5]{milnevalues} for the Brauer group and 
\cite{katri} for the Tate-Shafarevich group. 
\item The Tate-Shafarevich group $\Sh(\Pic^0_X)$ is finite if and only if 
$\Sh(\Alb_X)$ is finite \cite[I Remark 6.14 (c)]{adt}.
\end{enumerate}


{\it Proof of Theorem \ref{thm1}.} 
If Tate's conjecture holds for $X$ and the Tate-Shafarevich group of 
$\Alb_X$ is finite, then $(T\Br(\X))^*$ vanishes by Theorem \ref{maines}.
But the Brauer group is torsion with ${}_m\Br(\X)$ finite for 
every $m$, hence the vanishing of the Tate module implies that 
the $l$-primary part $\Br(\X)\{l\}$ is finite for every prime $l\not=p$. 

Conversely assume that the Brauer group $\Br(\X)$ is finite. Then 
$T\Br(\X)=0$, hence $(T\Sh(\Pic^0_X))^*=0$ by Theorem \ref{maines}. 
But the Tate-Shafarevich group is torsion with finite $m$-torsion 
for every $m$ \cite[I Remark 6.7]{adt}.
This implies that the $l$-primary part $\Sh(\Pic^0_X)\{l\}$ is finite 
for every prime $p\not=l$, hence that 
$\Sh(\Pic^0_X)$ is finite and 
then $\Sh(\Alb_X)$ are finite. Consequently $(\coker c)^*$ is finite
by Theorem \ref{maines}.
\endproof

{\it Proof of Proposition \ref{thm3}(1).}
First we consider $\ker \xi_d^{2d+1}$.
If $\Br(\X)$ is finite, then Tate's conjecture for divisors holds
on $\X$, hence Tate's conjecture for dimension one cycles 
holds on $\X$ \cite[Prop. 8.4]{milnevalues}, 
or equivalently $T_l H^{2d+1}_\et(\X,\Z(d))$ vanishes
for all $l\not=p$. This implies that $H^{2d+1}_\et(\X,\Z(d))\{l\}$
is the cotorsion of $H^{2d}_\et(\X,\Q_l/\Z_l(d))$, and this
is finite and vanishes for almost all $l$ by Gabber's theorem 
\cite{gabber}. We conclude that $H^{2d+1}_\et(\X,\ZZ(d))$ is finite,
and surjects onto $\ker l_d^{2d+1}$ in \eqref{mmdd}. By
\eqref{eins+} we see that $\coker \rho_d\to \ker \xi_d^{2d+1}$
has finite kernel and cokernel.

By Tate-Poitou duality $\coker \xi_d^{2d+1}$ is dual to
$H^2_\et(X^s,\hat \Z(1))^{G_K}$. Finiteness of $\Br(\X)$ implies
Tate's conjecture for $X$ in codimension one, hence this agrees
with the dual of $\NS(X)\otimes\hat\Z$, which is isomorphic to 
the dual of $\NS(X)$ as $\NS(X)$ is finitely generated. 
\endproof

In characteristic $0$, our method gives the following weaker result
(to obtain the full result one would need to prove 
Theorem \ref{localthm} in this situation).

\begin{theorem}
Let $\X$ be regular, proper and flat over the rings of integers of
a number field. If $\X$ has good reduction at all places above $p$
and if the $p$-primary component $\Br(\X)\{p\}$of the Brauer group 
is finite, then the $p$-primary component 
$\Sh(\Alb_X)\{p\} $ of the Tate-Shafarevich group of the albanese of 
the generic fiber is finite. 
\end{theorem}

\proof
If $\X$ has good reduction at $p$, then motivic cohomology agrees
with Sato's $p$-adic Tate twists, $\Z/p^r(n)\cong \T_r(n)$ 
by \cite{ichdede} and \cite[Thm. 4.8]{zhong}, see \cite[\S 1.4]{sato}. 
Hence the analog of Proposition \ref{cokerl}
can be proved by using \cite[Thm. 1.2.2]{sato}. On the other hand,
the vanishing of $H^{2d+1}_\et(\X_v,\Z_{(p)}(d))$ can be proved as in 
Corollary \ref{descc1} by using \cite[Thm. 1.3.1]{saisap}. 
Then the diagram chase in the diagram \eqref{eins+} gives the result.
\endproof

Finally, we show that the cokernel of $\rho_d$ is finite under the following 
conjecture, which was stated by Beilinson \cite[Conj. 5.2]{beilinson}
for number fields:

\begin{conjecture}\label{bbeeiill}
If $X^s$ smooth and projective over the separable closure $K^s$ of
a global field, then the albanese map 
$$\CH_0(X^s)^0 \to \Alb_{X^s}(K^s)$$
is an isomorphism.
\end{conjecture}

\begin{proposition}
Assuming Conjecture \ref{bbeeiill}, the cokernel of 
$$\rho_d: H^{2d}_\et(X,\ZZ(d))\to H^{2d}_\et(X^\sep,\ZZ(d))^{G_K}$$
is finite up to a $p$-group.
\end{proposition}

\proof
A norm argument show that the cokernel of $\rho_d$ is torsion, hence
it suffices to show that the target of $\rho_d$ is finitely generated. 
But by the Conjecture and Proposition \ref{rojtman}, 
it is an extension of $\ZZ$ by 
$(\CH_0(X^\sep)^0_{\Z'})^{G_K}\cong \Alb_{X}(K^\sep)_{\Z'}^{G_K}\cong
\Alb_X(K)_{\Z'}$, a finitely generated $\ZZ$-module.
\endproof

Finally we mention that our result implies the following Corollary,
of which the equivalence of 1)-4) is folklore.

\begin{corollary}
The first four statements are equivalent and imply the fifth statement:
\begin{enumerate}
\setlength\itemsep{0em}
\item Tate's conjecture holds for smooth and proper surfaces over finite fields. 
\item Tate's conjecture holds for divisors on smooth and projective varieties
over finite fields. 
\item Tate's conjecture holds for one-dimensional cycles on smooth and projective
varieties over finite fields. 
\item The Tate-Shafarevich group 
of any abelian variety $A$ over a global function field is finite.
\item Tate's conjecture holds
for divisors and for one-dimensional cycles on smooth and proper
varieties over global function fields for all $l\not=p$. 
\end{enumerate}
\end{corollary}

\proof
1) $\Leftrightarrow$ 3) follows by a Leftschetz theorem argument
\cite[Rem. 8.7]{milnevalues}, and 
1) $\Leftrightarrow$ 2) is proven in \cite[Thm. 4.3]{morrow}. 

The equivalence of 1) and 4) is known to the experts, 
but we repeat the argument for the convenience of the 
reader. Given a smooth and projective surface $\X$
over a finite field, it is explained in \cite[Proof of Thm. 1]{liu} 
how to obtain a surface satisfying the hypothesis of Theorem \ref{thm1}.
Then the finiteness of the Tate-Shafarevich group of the Jacobian
of the generic fiber implies the finiteness of $\Br(\X)$. 

Conversely, given an abelian variety $A$ over a function field $K$, we
can find another abelian variety $A'$ and a Jacobian $J$ of a smooth and proper
curve $X$ over $K$ which is isogeneous to $A\times A'$. Since 
$\Sh(J)$ and $\Sh(A)\oplus \Sh(A')$ differ by a finite group, it suffices 
by \cite{brauerIII} to observe that the Brauer group of a smooth and projective
model $\X$ of $X$ is finite by 1).

It remains to show that 2) implies 5). Let $X$ be a smooth and projective
variety over a function field $K$, and fix a prime $l$ different from $\chr K$.
By Gabber's theorem \ref{gabthm}, we can find a finite extension $K'$ 
of degree prime to $l$ such that
$X'=X\times_KK'$ has a model $\X$ which is regular and projective over the 
smooth and proper curve $C'$ with function field $K'$. 
Hence $\X$ is smooth and projective over a finite field,
and its Brauer group is finite by 2). From Theorem \ref{thm1}
we can then conclude that Tate's conjecture for divisors holds for $X'$,
which implies Tate's conjecture for $X$.
\endproof

\small

{\sc  Rikkyo University, Ikebukuro, Tokyo, Japan}

\textit{E-mail address:} {\tt geisser@rikkyo.ac.jp}

\begin{thebibliography}{SGA3}
\itemsep0em
\bibitem{beilinson} {\sc A.A.\ Beilinson}, Height pairings between algebraic
cycles. K-theory, arithmetic and geometry (Moscow, 1984--1986), 1--25, 
Lecture Notes in Math., 1289, Springer, Berlin, 1987. 


\bibitem{colthe} {\sc J.L. Colliot-Th\'el\`ene}, 
Cycles alg\'ebriques de torsion et K-th\'eorie alg\'ebrique.
Arithmetic algebraic geometry (Trento, 1991), 1--49, 
Lecture Notes in Math., 1553, Springer, Berlin, 1993. 

\bibitem{ct} {\sc J.L. Colliot-Th\'el\`ene}, Groupe de Chow des z\'ero-cycles 
sur les vari\'et\'es p-adiques (d'apr\`es S. Saito, K. Sato et al.). 
S\'eminaire Bourbaki. Vol. 2009/2010. Expos\'es 1012--1026. 
Ast\'erisque No. 339 (2011), Exp. No. 1012, vii, 1--30. 

\bibitem{dejong} {\sc J.\ A.\ de Jong}, Smoothness, semi-stability and alterations. 
Inst. Hautes \'Etudes Sci. Publ. Math. No. 83 (1996), 51--93. 


\bibitem{SGA4} {\sc P.\ Deligne}, La formule de dualite globale, in:
Th\'eorie des topos et cohomologie 
\'etale des sch\'emas. Tome 3. S\'eminaire de G\'eom\'etrie 
Alg\'ebrique du Bois-Marie 1963--1964. 
Dirig\'e par M. Artin, A. Grothendieck et J. L. Verdier. 
Lecture Notes in Mathematics, Vol. 305. Springer-Verlag, Berlin-New York, 1973. vi+640 pp.


\bibitem{gabber} {\sc O.\ Gabber}, Sur la torsion dans la cohomologie 
$l$-adique d'une vari\'et\`e. 
C. R. Acad. Sci. Paris S\'er. I Math. 297 (1983), no. 3, 179--182. 



\bibitem{ichdede} {\sc T.\ Geisser}, 
Motivic cohomology over Dedekind rings. Math. Z. 248 (2004), no. 4, 773--794. 


\bibitem{ichduality} {\sc T.\ Geisser}, Duality via cycle complexes,
Ann.\  of Math.\  (2) 172 (2010), no.\  2, 1095--1126.

\bibitem{ichstructure} {\sc T.\ Geisser}, On the structure of 
etale motivic cohomology. J. Pure Appl. Algebra 221 (2017), no. 7, 1614--1628. 

\bibitem{ichtata} {\sc T.\ Geisser}, Duality of integral  \'etale
motivic cohomology, to appear in Proc. International Colloquium 
in K-theory TIFR in January, 2016. https://arxiv.org/abs/1712.09021


\bibitem{ichsurface} {\sc T.\ Geisser}, Comparing the Brauer group
and the Tate-Shafarevich group, to appear in:
J. Inst. Math. Jussieu.


\bibitem{ichmarcII} {\sc T.\ Geisser, M.\ Levine},
The Bloch-Kato conjecture and a theorem of Suslin-Voevodsky, 
J. reine angew. Math. 530 (2001), 55--103.



\bibitem{tan} {\sc C.D.\ Gonz\'alez-Avil\'es, Ki-Seng Tan}, 
A generalization of the Cassels-Tate dual exact sequence. 
Math. Res. Lett. 14 (2007), no. 2, 295--302.

\bibitem{brauerIII} {\sc A.\ Grothendieck}, Le groupe de Brauer. III. 
Exemples et compl\'ements. Dix expos\'es sur la cohomologie des sch\'emas, 
88--188, Adv. Stud. Pure Math., 3, North-Holland, Amsterdam, 1968.

\bibitem{gabberthm} {\sc L.\ Illusie, Y.\ Laszlo, F.\ Orgogozo}. 
Travaux de Gabber sur l'uniformisation locale et la cohomologie  \'etale des 
sch\'emas quasi-excellents.
S\'eminaire \`a l'\'Ecole Polytechnique 2006--2008. 
Ast\'erisque No. 363--364 (2014). Soci\'et\'e Math\'ematique de 
France, Paris, 2014. 

\bibitem{jannsenggq} {\sc U.\ Jannsen}, On the $l$-adic cohomology of 
varieties over number fields and its Galois cohomology. 
Galois groups over $\Q$ (Berkeley, CA, 1987), 315--360, 
Math. Sci. Res. Inst. Publ., 16, Springer, New York, 1989.


\bibitem{jss} {\sc U.\ Jannsen, S.\ Saito, K.\ Sato}, 
\'Etale duality for constructible sheaves on arithmetic schemes.
J. Reine Angew. Math. 688 (2014), 1--65. 


\bibitem{kahnhandbook} {\sc B.\ Kahn}, 
Algebraic K-theory, algebraic cycles and arithmetic geometry. 
Handbook of K-theory. Vol. 1, 2, 351--428, Springer, Berlin, 2005.




\bibitem{katri}{\sc K.\ Kato, F. Trihan}, On the conjectures of Birch and
Swinnerton-Dyer in characteristic $p> 0$.
Invent. Math. 153 (2003), no. 3, 537--592.

\bibitem{KS} { M.\ Kerz, S.\ Saito}, 
Cohomological Hasse principle and motivic cohomology
for arithmetic schemes. 
Publ.\  Math.\  Inst.\  Hautes Etudes Sci.\  115 (2012), 123--183.



\bibitem{levine} {\sc M.\ Levine}, K-theory and motivic cohomology of schemes,
https://faculty.math.illinois.edu/K-theory/0336/mot.pdf.

\bibitem{lichtenbaumcurve} {\sc S.\ Lichtenbaum}, Duality theorems for curves 
over p-adic fields. Invent. Math. 7 (1969), 120--136. 

\bibitem{liu}{\sc Q. Liu, D. Lorenzini, M. Raynaud},
On the Brauer group of a surface. Invent. Math. 159 (2005), no. 3, 673--676.


\bibitem{milnevalues} {\sc J.\ Milne}, 
Values of zeta functions of varieties over finite fields. 
Amer. J. Math. 108 (1986), no. 2, 297--360.

\bibitem{adt} {\sc J.S.\ Milne}, Arithmetic duality theorems. Second edition. 
BookSurge, LLC, Charleston, SC, 2006. viii+339 pp. ISBN: 1-4196-4274-X 


\bibitem{morrow} {\sc M.\ Morrow}, 
A Variational Tate Conjecture in crystalline cohomology, 
https://arxiv.org/pdf/1408.6783.pdf 


\bibitem{saisa} {\sc S.\ Saito, K.\ Sato}, A finiteness theorem for zero-cycles 
over $p$-adic fields. With an appendix by Uwe Jannsen. 
Ann. of Math. (2) 172 (2010), no. 3, 1593--1639.

\bibitem{saisap} {\sc S.\ Saito, K.\ Sato},
Zero-cycles on varieties over p-adic fields and Brauer groups. 
Ann. Sci. \'Ec. Norm. Sup\'er. (4) 47 (2014), no. 3, 505--537.

\bibitem{sato} {\sc K.\ Sato}, $p$-adic \'etale Tate twists and arithmetic duality.
With an appendix by Kei Hagihara. Ann. Sci. \'Ecole Norm. Sup. (4) 40 (2007), 
no. 4, 519--588.



\bibitem{thomason} {\sc R.W.\ Thomason}, 
Algebraic K-theory and \'etale cohomology. 
Ann. Sci. \'Ecole Norm. Sup. (4) 18 (1985), no. 3, 437--552. 


\bibitem{zhong}  {\sc C.\ Zhong}, 
Comparison of dualizing complexes. J. Reine Angew. Math. 695 (2014), 1--39.
\end{thebibliography}
\end{document}